\documentclass[12pt]{article}
\usepackage{graphicx,float,subfig}               
\usepackage{amsmath,pgf,tikz,tkz-graph}
\usetikzlibrary{arrows}
\usepackage{graphicx}
\usepackage{amssymb}
\usepackage{epstopdf}
\usepackage[noblocks]{authblk}
\newtheorem{theorem}{Theorem}[section]

\title{Some New Methods for Constructing 4-critical Planar Graphs}
\author{Zhou Guofei\thanks{gfzhou@nju.edu.cn}\\
Department of Nanjing University, Nanjing 210093, China}

\date{}                                         
\begin{document}
\baselineskip 0.29in
\maketitle
\begin{abstract}
 A graph $G$ is said to be $k$-critical if $G$ is $k$-colorable and $G-e$ is not $k$-colorable for every edge 
$e$ of $G$. In this paper, we present some new methods from two or more small 4-critical graphs to  construct a larger 4-critical planar graphs. 
\end{abstract}
\section{Introduction} 

Let $G=(V,E)$ be a graph, $G$ is said to be $k$-colorable if there is a assignment of $k$ colors to the vertices of $G$ such that no two adjacent vertices of $G$ get the same color.  The chromatic number of $G$, denoted by $\chi(G)$, is the least integer $k$ such that $G$ is $k$-colorable.
A graph $G$ is said to be $k$-critical if $G$ is $k$-colorable and $G-e$ is not $k$-colorable for every edge 
$e$ of $G$. A planar graph is a graph that can be embedded in the plane, i.e., it can be drawn on the plane in such a way that its edges intersect only at their endpoints.
About 1950 G. A. Dirac introduced the concept of criticality as a methodological means in the theory of graph colouring. He himself, G. Hajbs and T. Gallai were the first to develop special constructions for creating colour-critical graphs and establishing theorems on their properties. 
The basic tool to be used to construct critical graphs is Haj\'{o}s’ construction \cite{hajos, TB}: Let $G_1$ and $G_2$ be two disjoint graphs with 
edges $x_1y_1$ and $x_2y_2$, let $G$ be a new graphs obtained from $G_1$ and $G_2$ by removing $x_1y_1$ and $x_2y_2$, identifying $x_1$ and $x_2$, and joining $y_1$ and $y_2$
by a new edge. It is well known that if $G_1$ and $G_2$ are $k$-critical graphs, so is $G$.

There is a survey on methods of constructing critical graphs, see \cite{SS}. As for 4-critical planar graphs, 
there are only a few methods which are known to constructing 4-critical planar graphs 
(see \cite{B.Zhou,AZ,Bor,BG1,Koester,Koester1} etc). In this paper, we present some new methods 
for constructing 4-critical planar graphs, given two or more 4-critical planar graphs. 

\section{Main Results}
{\bf  Some new methods for constructing 4-critical planar graphs}

Let $G_1, G_2, G_3$ be three planar graphs. Suppose that $u_1u_2u_3$, $v_1v_2v_3$ and $w_1w_2w_3$ are three paths each lie in the boundary of some face of $G_1$, $G_2$ and $G_3$ respectively.  Let $G=G_1\circ G_2\circ G_3$ be a graph (see Figure \ref{fig1}) obtained by

(a) Deleting the edges $u_2u_3$, $v_2v_3$ and $w_2w_3$;

(b) Identifying  $u_1$ with $w_3$, $u_3$ with $v_1$ and $v_1$ with $w_1$;

(c) Adding a new vertex $z$, and adding three new edges $zu_2$, $zv_2$, $zw_2$.

\begin{figure}[H]
\centering
\subfloat[]{
\begin{tikzpicture}[scale=1]
\renewcommand*{\VertexSmallMinSize}{1.pt}
\renewcommand*{\VertexInnerSep}{1.pt}
\GraphInit[vstyle=Classic]
\Vertex[x=0.57, y=3.48,Lpos=80]{$w_3$}
\Vertex[x=0.43, y=3.48,Lpos=100]{$u_1$}
\Vertex[x=-.25, y=2.3]{$u_2$}
\Vertex[x=-0.93, y=1.12,Lpos=180]{$u_3$}
\Vertex[x=-0.86, y=1,Lpos=-90]{$v_1$}
\Vertex[x=.5, y=1,Lpos=90]{$v_2$}
\Vertex[x=1.86, y=1,Lpos=-88]{$v_3$}
\Vertex[x=1.93, y=1.12]{$w_1$}
\Vertex[x=1.25, y=2.3,Lpos=180]{$w_2$}
\Edges($u_1$,$u_2$,$u_3$)
\Edges($w_1$,$w_2$,$w_3$)
\Edges($v_1$,$v_2$,$v_3$)
\draw (0.57,3.48) .. controls (2.99,3.31) .. (1.93,1.12);
\draw  (0.43,3.48) .. controls (-1.97,3.31) .. (-0.93,1.12);
\draw  (-0.86,1) .. controls (0.5,-0.76) .. (1.86,1);
\draw (-1.11,2.8) node{$G_1$};
\draw (0.5,0.12) node{$G_2$};
\draw (2.12,2.8) node{$G_3$};
\end{tikzpicture}
}
\subfloat[]{
\begin{tikzpicture}[scale=1]
\renewcommand*{\VertexSmallMinSize}{1.2pt}
\renewcommand*{\VertexInnerSep}{1.2pt}
\GraphInit[vstyle=Classic]
\Vertex[x=0.5, y=3.6,Lpos=90,Math,L=u_1]{u_1}
\Vertex[x=-1, y=1,Lpos=180,Math]{v_1}
\Vertex[x=2, y=1,Math]{w_1}
\Vertex[x=.5, y=1,Lpos=135,Math]{v_2}
\Vertex[x=1.25, y=2.3,Lpos=180,Math]{w_2}
\Vertex[x=-.25, y=2.3,Math]{u_2}
\Vertex[x=.5, y=1.87,Math]{z}
\AddVertexColor{white}{z}
\Edges(u_1,u_2)
\Edges(w_1,w_2)
\Edges(v_1,v_2)
\Edges(z,v_2)
\Edges(z,w_2)
\Edges(z,u_2)
\draw (u_1) .. controls (2.99,3.31) .. (w_1);
\draw (u_1) .. controls (-1.97,3.31) .. (v_1);
\draw (v_1) .. controls (0.5,-0.76) .. (w_1);
\draw (-1.11,2.8) node{$G_1$};
\draw (0.5,0.12) node{$G_2$};
\draw (2.12,2.8) node{$G_3$};

\draw [dash pattern=on 3pt off 3pt] (u_1)-- (w_2);
\draw [dash pattern=on 3pt off 3pt] (w_1)-- (v_2);
\draw [dash pattern=on 3pt off 3pt] (v_1)-- (u_2);
\end{tikzpicture}
}
\caption{$G=G_1\circ G_2\circ G_{3}$}
\label{fig1}
\end{figure}
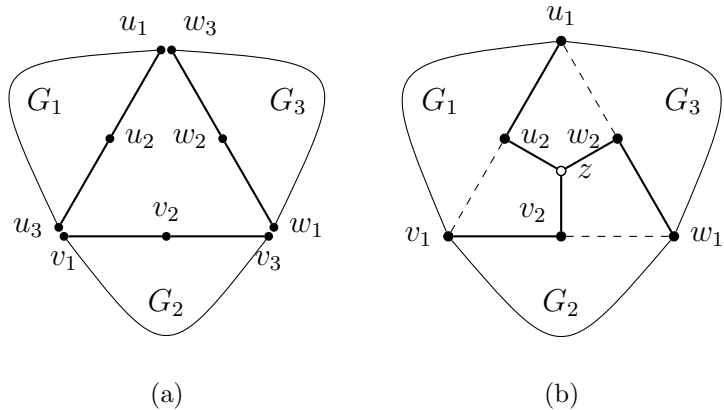

\begin{figure}[H]
\centering
\subfloat{
\begin{tikzpicture}[scale=0.6]
\renewcommand*{\VertexSmallMinSize}{1.pt}
\renewcommand*{\VertexInnerSep}{1.pt}
\GraphInit[vstyle=Classic]
\Vertex[x=8.02, y=4.18,Lpos=90,Math,L=u_1]{u_1}
\Vertex[x=5, y=2,Lpos=180,Math]{v_1}
\Vertex[x=6.51, y=3.09,Lpos=0,Math]{u_2}
\Vertex[x=6.14, y=-1.54,Lpos=-90,Math]{w_1}
\Vertex[x=9.86, y=-1.55,Lpos=-45,,Ldist=-3pt,Math]{x_1}
\Vertex[x=11.04, y=1.97,Math]{y_1}
\Vertex[x=5.57, y=0.23,Lpos=87.5,Ldist=-1pt,Math]{v_2}
\Vertex[x=8, y=-1.55,Lpos=45,Ldist=-3pt,Math]{w_2}
\Vertex[x=10.45, y=0.21,Lpos=95,Math]{x_2}
\Vertex[x=9.52, y=3.08,Lpos=-90,Math]{y_2}
\Vertex[x=8.01, y=1.01,Ldist=-1pt,Lpos=90,Math]{z}
\Edges(u_2,v_1)
\Edges(v_2,w_1)
\Edges(w_2,x_1)
\Edges(x_2,y_1)
\Edges(y_2,u_1)
\Edges(z,u_2)
\Edges(z,v_2)
\Edges(z,w_2)
\Edges(z,x_2)
\Edges(z,y_2)
\AddVertexColor{white}{z}
\draw [dash pattern=on 3pt off 3pt] (u_1)-- (u_2);
\draw [dash pattern=on 3pt off 3pt] (v_1)-- (v_2);
\draw [dash pattern=on 3pt off 3pt] (w_1)-- (w_2);
\draw [dash pattern=on 3pt off 3pt] (x_1)-- (x_2);
\draw [dash pattern=on 3pt off 3pt] (y_1)-- (y_2);
\draw (u_1) .. controls (11,5) .. (y_1);
\draw (u_1) .. controls (5,5) .. (v_1);
\draw (v_1) .. controls (3,-1) .. (w_1);
\draw (w_1) .. controls (8,-4) .. (x_1);
\draw (x_1) .. controls (13,-1) .. (y_1);
\draw (6,4) node{$G_1$};
\draw (4.44,-0.36) node{$G_2$};
\draw (8,-3) node{$G_3$};
\draw (11.84,-0.32) node{$G_4$};
\draw (10,4) node{$G_5$};
\end{tikzpicture}}
\caption{The graph $G=G_1\circ G_2\circ G_{3}\circ G_4\circ G_5$}
\label{fig2}
\end{figure}
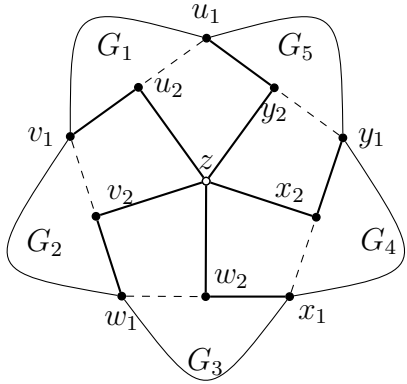

\begin{figure}[H]
\centering
\includegraphics[scale=1.7]{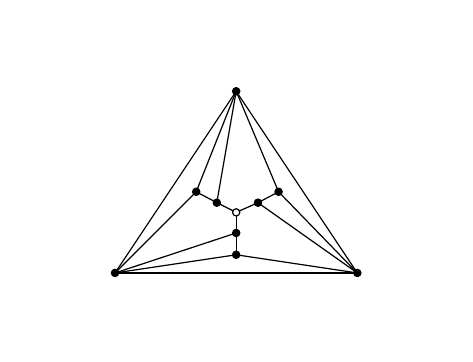}
\caption{The graph $G=K_4\circ K_4\circ K_4$}
\label{example1}
\end{figure}

If $G_1=G_2=G_3=K_4$, the new graph $G=G_1\circ G_2\circ G_3$ is shown in Figure \ref{example1}. It is a 4-critical planar graph.

\begin{figure}[H]
\centering
\includegraphics[scale=1.]{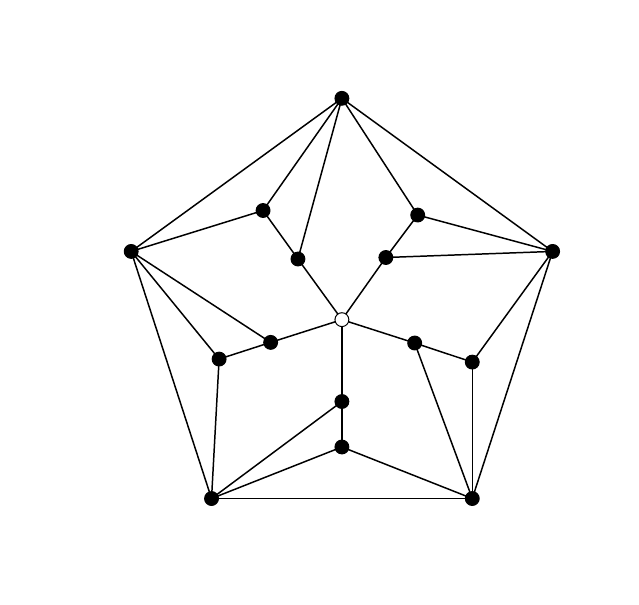}
\caption{The graph $G=K_4\circ K_4\circ K_4\circ K_4\circ K_4$}
\label{example2}
\end{figure}

If $G_1=G_2=G_3=G_4=G_5=K_4$, the new graph $G=G_1\circ G_2\circ G_3\circ G_4\circ G_5$ is shown in Figure \ref{example2}. It is a 4-critical planar graph.

\begin{theorem}\label{thm1}
If $G_1$, $G_2$ and $G_3$ are 4-critical planar graphs, then $G=G_1\circ G_2\circ G_3$ is a 4-critical graph.
\end{theorem}

{\bf Proof.} (i) It is clear that $G$ is still a planar graph, by the Four Color Theorem, we have $\chi (G)\le 4$. We next show that $\chi (G)\ge 4$.

If $G$ can be colored with 3 colors,  let $c$ be a 3-coloring of $G$. Since $G_1$ is 4-critical, $G_1-w_2w_3$ can be colored with 3 colors. In such a 
coloring, $w_3$ and $w_2$ must be colored with the same color (otherwise this would be a 3-coloring of $G_1$
as well). This implies that $c(u_1)=c(w_2)$. Similarly, we have that $c(w_1)=c(v_2)$ and $c(v_1)=c(u_2)$. 
Note that $z$ is adjacent to each vertex of $\{u_2,v_2,w_2\}$, at least two vertices of them have the same color, 
without loss of generality, we assume that $c(u_2)=c(v_2)$, 
then $c(v_1)=c(v_2)$ since $c(u_2)=c(v_1)$, this is impossible since $v_1$ and $v_2$ are adjacent in $G$. 
Therefore, we have $\chi (G)= 4$. 

\medskip
(ii) We show that $\chi (G-e)\le 3$ for every edge in $G$. Let $H_1,H_2$ and $H_3$ be the subgraph of $G$ induced by $V(G_1), V(G_2)$ and $V(G_3)$ respectively. It is clear that $H_1=G_1-u_2u_3, H_2=G_2-v_2v_3$ and $H_3=G_3-w_2w_3$ respectively. Consider the following three cases:

\noindent {\bf Case 1.} $e\in\{zu_2,zv_2,zw_2\}$.

Since $G_1$ is 4-critical, there is a three coloring $c_1$ of $G_1-{u_2u_3}$ such that $u_2$ and $u_3$ get the same 
color. Without loss of generality, assume that $c_1(u_2)=c_1(u_3)=1$ and $c_1(u_1)=2$. Also, since $G_2$ is 4-critical, 
there is a three coloring $c_2$ of $G_2-{v_2v_3}$ such that $c_1(v_1)=1$ and $c_2(v_2)=c_1(v_3)=3$. (if $c_2(v_2)=c_2(v_3)=2$, we can swap the colors of 2 and 3 in $G_2$). Similarly, since $G_3$ is 4-critical, 
there is a three coloring $c_3$ of $G_3-{w_2w_3}$ such that $c_3(w_1)=3$ and $c_3(w_2)=c_3(w_3)=2$. As $u_1=w_3$, $v_1=u_3$ and $w_1=v_3$, we see that $c_1,c_2,c_3$ form a 3-coloring of $G-z$, this coloring can
obviously be extended to a 3-coloring of $G-e$.

\noindent {\bf Case 2.} $e\in E(G)/\{zu_2,zv_2,zw_2,u_1u_2,v_1v_2,w_1w_2\}$.

Without loss of generality, assume that $e\in E(G_1)-u_1u_2$.  Since $G_1$ is 4-critical, there is a 3-coloring  of $G_1-e$. This coloring induces a 3-coloring $c_1$ of $H_1-e$. Without loss of generality, we assume that $c_1(v_1)=1$ and $c_1(u_2)=2$. 

If $c_1(u_1)=3$, then since there is a 3-coloring $c_3$ of $H_3$ such that $c_3(u_1)=c_3(w_2)=3$ and $c_3(w_1)=2$. Similarly, there is a 3-coloring $c_2$ of $H_2$ such that $c_2(w_1)=c_2(v_2)=2$ and $c_2(v_1)=1$. Therefore, $c_1,c_2,c_3$ induce a 3-coloring of $G-z$. This coloring can be extended to a 3-coloring of $G-e$, since there are two
colors used in the neighbors of $z$.

If $c_1(u_1)=3$, by a similar argument as above, we can show that $G-e$ is 3-colorable. 

\noindent {\bf Case 3.} $e\in\{u_1u_2,v_1v_2,w_1w_2\}$.

Without loss of generality, assume that $e=u_1u_2$. Since $G_1$ is 4-critical, there is a 3-coloring of $G_1-e$ such that $u_1$ and $u_2$ have the same color. This coloring induces a 3-coloring $c_1$ of $H_1-e$. Assume, without loss of generality, $c_1(u_1)=c_1(u_2)=1$ and $c_1(v_1)=2$. Since $G_3$ is 4-critical, there is a 3-coloring $c_3$ of $H_3$ such that $c_3(u_1)=c_3(w_2)=1$ and $c_3(w_1)=3$. Similarly, there is a 3-coloring $c_2$ of $H_2$ such that $c_2(w_1)=c_2(v_2)=3$ and $c_2(v_1)=2$. Therefore, $c_1,c_2,c_3$ induce a 3-coloring of $G-z$. This coloring can be extended to a 3-coloring of $G-e$, since there are two colors used in the neighbors of $z$. \hfill $\Box$

\medskip
Let $G_1, G_2,\cdots G_{2k+1}$ be three planar graphs. Suppose that $u^1_1u^1_2u^1_3$, 
$u^2_1u^2_2u^2_3$, $\cdots$,
$u^{2k+1}_1u^{2k+1}_2u^{2k+1}_3$  are $(2k+1)$ paths each lie in the boundary of some face of $G_1$, $G_2$ and $G_{2k+1}$ respectively.  Let $G=G_1\circ G_2\circ\cdots\circ G_{2k+1}$ be a graph (if $k=2$, see Figure \ref{fig2}) obtained by

(a) Deleting the edges $u^1_2u^1_3$, $u^2_2u^2_3$,$\cdots$, $u^{2k+1}_2u^{2k+1}_3$;

(b) Identifying  $u^1_1$ with $u^{2k+1}_3$, $u^1_3$ with $u^2_1$ and $u^2_3$ with $u^3_1$, $\cdots$,$u^{2k}_3$ with $u^{2k+1}_1$;

(c) Adding a new vertex $z$, and adding $2k+1$ new edges $zu^1_2$, $zu^2_2$, $\cdots$, $zu^{2k+1}_2$.

\begin{theorem}\label{thm2}
If $G_1$, $G_2$ and $G_{2k+1}$ are 4-critical planar graphs, then $G=G_1\circ G_2\circ\cdots\circ G_{2k+1}$ is a 4-critical planar graph.
\end{theorem}

{\bf Proof.} (i) It is clear that $G$ is still a planar graph, by the Four Color Theorem, we have $\chi (G)\le 4$. We next show that $\chi (G)\ge 4$.

Since each $G_i$ is 4-critical, $G_i-u^i_2u^i_3$ ($1\le i\le 2k+1$) is 3-colorable and $c(u^i_2)=c(u^i_3)$ for any 3-coloring of $G_i-u^i_2u^i_3$. If $G$ is 3-colorable, without loss of generality, let $c$ be a 3-coloring of $G_i-u^i_2u^i_3$.
Denote by $H$ the subgraph of $G$ induced by the vertex set $\{z,u^1_1,u^1_2,u^1_3,u^2_2,u^2_3,\cdots,u^{2k}_{2},u^{2k}_3,u^{2k+1}_2\}$. Then $H$ is 3-colorable, but this is impossible since $c(u^i_2)=c(u^i_3)$ ($1\le i\le 2k+1$).

\medskip
(ii) We show that $\chi (G-e)\le 3$ for every edge in $G$. Let $H_1,H_2$,$\cdots$, $H_{2k+1}$ be the subgraph of $G$ induced by $V(G_1), V(G_2)$,$\cdots$, $V(G_{2k+1})$ respectively. It is clear that $H_i=G_i-u^i_2u^i_3$ ($1\le i\le 2k+1$). Consider the following three cases:

\noindent {\bf Case 1.} $e\in\{zu^1_2,zu^2_2,\cdots,zu^{2k+1}_2\}$.

Without loss of generality, assume that $e=zu^1_2$. 
Since $G_1$ is 4-critical, there is a three coloring $c_1$ of $G_1-{u^1_2u^1_3}$ such that $c_1(u^1_2)=c_1(u^1_3)=1$ and $c_1(u^1_1)=3$. Also, since $G_2$ is 4-critical, 
there is a three coloring $c_2$ of $G_2-{u^2_2u^2_3}$ such that $c_2(u^1_3)=1$ and $c_2(u^2_2)=c_2(u^2_3)=2$. 
In general, if $i\ge 3$ is odd, since $G_i$ is 4-critical, 
there is a three coloring $c_i$ of $G_i-{u^i_2u^i_3}$ such that $c_i(u^{i-1}_3)=2$ and $c_i(u^i_2)=c_i(u^i_3)=3$;
if $i\ge 3$ is even, since $G_i$ is 4-critical, 
there is a three coloring $c_i$ of $G_i-{u^i_2u^i_3}$ such that $c_i(u^{i-1}_3)=3$ and $c_i(u^i_2)=c_i(u^i_3)=2$;
Finally, we color $z$ with color 1. We see that $c_1,c_2,\cdots,c_{2k+1}$ together with the color of $z$ extends  a 3-coloring of $G-e$.

\noindent {\bf Case 2.} $e\in E(H_i)$, for ($1\le i\le 2k+1$).

Without loss of generality, assume that $e\in E(H_1)$.  

\medskip

Suppose first that $e=u^1_1u^1_2$. Since $G_1$ is 4-critical, there is a 3-coloring of $G_1-e$ such that $u^1_1$ and $u^1_2$ have the same color. This coloring induces a 3-coloring $c_1$ of $H_1-e$. Assume, without loss of generality, $c_1(u^1_1)=c_1(u^1_2)=1$ and $c_1(u^1_3)=3$. Since $G_2$ is 4-critical, there is a three coloring $c_i$ of $G_2-{u^2_2u^2_3}$ such that $c_2(u^{1}_3)=3$ and $c_2(u^2_2)=c_2(u^2_3)=2$; In general, if $i\ge 3$ is odd, since $G_i$ is 4-critical, 
there is a three coloring $c_i$ of $G_i-{u^i_2u^i_3}$ such that $c_i(u^{i-1}_3)=2$ and $c_i(u^i_2)=c_i(u^i_3)=1$;
if $i\ge 3$ is even, since $G_i$ is 4-critical, 
there is a three coloring $c_i$ of $G_i-{u^i_2u^i_3}$ such that $c_i(u^{i-1}_3)=1$ and $c_i(u^i_2)=c_i(u^i_3)=2$;
Finally, we color $z$ with color 3. We see that $c_1,c_2,\cdots,c_{2k+1}$ together with the color of $z$ extends  a 3-coloring of $G-e$.
\medskip

Suppose now that $e\in H_1-u^1_1u^1_2$. 

Since $G_1$ is 4-critical, there is a 3-coloring  of $G_1-e$. This coloring induces a 3-coloring $c_1$ of $H_1-e$. 
Assume without loss of generality, let $c_1(u^1_1)=2$, $c_1(u^1_2)=1$ and $c_1(u^1_3)=k$, where $k\in\{2,3\}$.

since $G_2$ is 4-critical, 
there is a three coloring $c_2$ of $G_2-{u^2_2u^2_3}$ such that $c_2(u^1_3)=k$ and $c_2(u^2_2)=c_2(u^2_3)=1$;

If $i\ge 3$ is odd, since $G_i$ is 4-critical, 
there is a three coloring $c_i$ of $G_i-{u^i_2u^i_3}$ such that $c_i(u^{i-1}_3)=1$ and $c_i(u^i_2)=c_i(u^i_3)=2$;
if $i\ge 3$ is even, since $G_i$ is 4-critical, 
there is a three coloring $c_i$ of $G_i-{u^i_2u^i_3}$ such that $c_i(u^{i-1}_3)=2$ and $c_i(u^i_2)=c_i(u^i_3)=1$;
Finally, we color $z$ with color 3. We see that $c_1,c_2,\cdots,c_{2k+1}$ together with the color of $z$ extends  a 3-coloring of $G-e$. \hfill $\Box$

\begin{figure}[H]
\centering
\subfloat{
\begin{tikzpicture}[scale=1]
\renewcommand*{\VertexSmallMinSize}{1.5pt}
\renewcommand*{\VertexInnerSep}{1.5pt}
\GraphInit[vstyle=Classic]
\Vertex[x=1, y=5,Lpos=90,Math,L=u_1]{A}
\Vertex[x=-1.85, y=0.01,Lpos=180,Math,L=u_3]{B}
\Vertex[x=3.87, y=0.01,Lpos=0,Math,L=v_3]{C}
\Vertex[x=-0.43, y=2.5,Lpos=0,Math,L=u_2]{D}
\Vertex[x=2.44, y=2.5,Lpos=180,Math,L=v_2]{E}
\tikzset{VertexStyle/.append  style={fill=white}}
\Vertex[x=1, y=1.68,Lpos=-90,Math,L=z]{H}
\Edges(A,E)
\Edges(C,B)
\Edges(D,A)
\Edges(D,H,E)
\Edge(H)(A)
\tikzset{EdgeStyle/.style = {-,bend left=80,color=gray}}
\Edge(A)(C)
\tikzset{EdgeStyle/.style = {-,bend right=80,color=gray}}
\Edge(A)(B)
\draw [dash pattern=on 3pt off 3pt] (E)-- (C);
\draw [dash pattern=on 3pt off 3pt] (B)-- (D);
\draw (3.5,2.5) node{$G_2$};
\draw (-1.5,2.5) node{$G_1$};
\end{tikzpicture}}
\caption{The graph $G=g_3(G_1,G_2)$}
\label{fig3}
\end{figure}
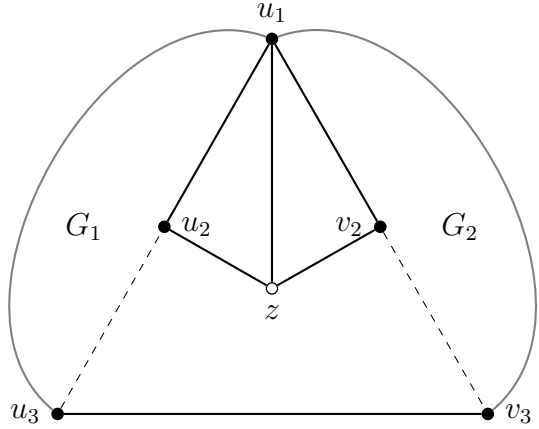

Let $G_1, G_2$ be two planar graphs. Suppose that $u_1u_2u_3$ and $v_1v_2v_3$ are two paths each lie in the boundary of some face of $G_1$, $G_2$ respectively.  Let $G=g_3(G_1, G_2)$ be a graph (see Figure \ref{fig3}) obtained from $G_1$ and $G_2$ by

(a) Deleting the edges $u_2u_3$ and $v_2v_3$;

(b) Identifying  $u_1$ with $v_1$;

(c) Adding a new vertex $z$, and adding three new edges $zu_1$, $zu_2$, $zv_2$;

(d) Adding a new edge $u_3v_3$.

\begin{figure}[H]
\centering
\includegraphics[scale=1.]{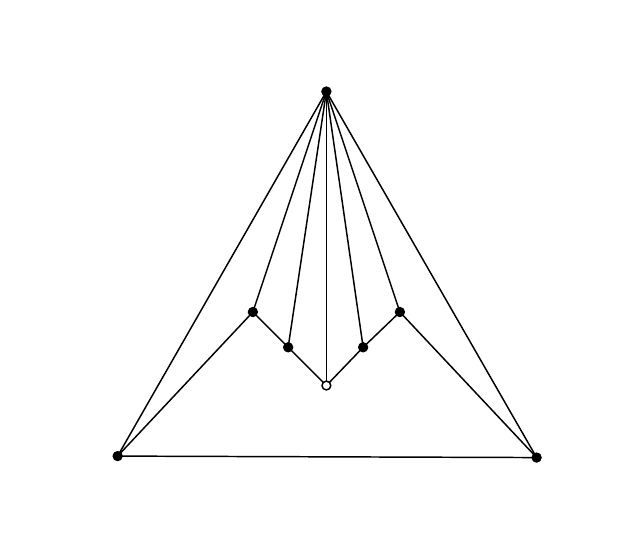}
\caption{The graph $G=g_3(K_4, K_4)$}
\label{example3}
\end{figure}

\medskip 
If $G_1=G_2=K_4$, for example, $G=g_3(K_4,K_4)$, see Figure \ref{example3}, it is a 4-critical planar graph.

\begin{theorem}\label{thm2}
If $G_1$, $G_2$ and $G_{2k+1}$ are 4-critical planar graphs, then $G=g_3(G_1, G_2)$ is a 4-critical planar graph.
\end{theorem}

{\bf Proof.} (i) It is clear that $G$ is still a planar graph, by the Four Color Theorem, we have $\chi (G)\le 4$. We next show that $\chi (G)\ge 4$.

If $G$ is 3-colorable, then
since $G_1$ and $G_2$ are 4-critical, $G_1-u_2u_3$ and $G-v_2v_3$ are 3-colorable and $c(u_2)=c(u_3)$ and  $c(v_2)=c(v_3)$ for any 3-coloring $c$ of $G$.  Note that $u_3$ and $v_3$ are adjacent, without loss of generality, let $c(u_1)=1$, $c(u_2)=c(u_3)=2$ and $c(v_2)=c(v_3)=3$. This is impossible since we can not color $z$ with color 1,2 or 3. 

(ii) We show that $\chi (G-e)\le 3$ for every edge in $G$. Let $H_1,H_2$  be the subgraph of $G$ induced by $V(G_1), V(G_2)$  respectively. It is clear that $H_1=G_1-u_2u_3, H_2=G_2-v_2v_3$. Consider the following  cases:

\medskip

{\bf Case 1.} $e\not\in (E(H_1)\cup E(H_2))$.

If $e=u_3v_3$, let $c_1$ and $c_2$ be the two 3-colorings of $H_1$ and $H_2$ respectively. Since $c_1(u_2)=c_1(u_3)$ and $c_2(v_2)=c_2(v_3)$, without loss of generality, let $c_1(u_1)=c_2(u_1)=1$, $c_1(u_2)=c_1(u_3)=2$ and $c_2(v_2)=c_2(v_3)=2$, then we can color $z$ with color 3. Therefore, $c_1,c_2$
together with the color of $z$ extend to a 3-coloring of $G$. 

If $e=zu_1$, let $c_1$ and $c_2$ be the two 3-colorings of $H_1$ and $H_2$ respectively. Since $c_1(u_2)=c_1(u_3)$ and $c_2(v_2)=c_2(v_3)$, without loss of generality, let $c_1(u_1)=c_2(u_1)=1$, $c_1(u_2)=c_1(u_3)=2$ and $c_2(v_2)=c_2(v_3)=3$, then we can color $z$ with color 1. Therefore, $c_1,c_2$
together with the color of $z$ extend to a 3-coloring of $G$. 

If $e=zu_2$ or $e=zv_2$, without loss of generality, assume that $e=zu_2$. Let $c_1$ and $c_2$ be the two 3-colorings of $H_1$ and $H_2$ respectively. Since $c_1(u_2)=c_1(u_3)$ and $c_2(v_2)=c_2(v_3)$, without loss of generality, let $c_1(u_1)=c_2(u_1)=1$, $c_1(u_2)=c_1(u_3)=2$ and $c_2(v_2)=c_2(v_3)=3$, then we can color $z$ with color 2. Therefore, $c_1,c_2$ together with the color of $z$ extend to a 3-coloring of $G$. 

\medskip

{\bf Case 2.} $e\in (E(H_1)\cup E(H_2))$.

Without loss of generality, assume that $e\in E(H_1)$. 

If $e=u_1u_2$, let $c_2$ be the  3-colorings  $H_2$. Since $G_1$ is 4-critical, there is a 3-coloring $c_1$ of $G_1-u_1u_2$ such that $c_1(u_1)=c_1(u_2)$. Without loss of generality, let $c_1(u_1)=c_1(u_2)=c_2(u_1)=1$, $c_1(u_3)=2$ and $c_2(v_2)=c_2(v_3)=3$, then we can color $z$ with color 2. Therefore, $c_1,c_2$ together with the color of $z$ extend to a 3-coloring of $G$. 

If $e\not=u_1u_2$, let $c_2$ be the  3-colorings  $H_2$. Since $G_1$ is 4-critical, there is a 3-coloring $c_1$ of $G_1-e$. This implies that $c_1(u_2)\ne c_1(u_3)$. Without loss of generality, let $c_1(u_1)=c_2(u_1)=1$, $c_1(u_3)=k$ ($k\in\{1,3\}$) and $c_2(v_2)=c_2(v_3)=2$, then we can color $z$ with color 3. Therefore, $c_1,c_2$ together with the color of $z$ extend to a 3-coloring of $G$. \hfill $\Box$

\begin{figure}[H]
\centering
\subfloat[]{
\begin{tikzpicture}[x=2cm,y=2cm]
\renewcommand*{\VertexSmallMinSize}{1.2pt}
\renewcommand*{\VertexInnerSep}{1.2pt}
\GraphInit[vstyle=Classic]
\Vertex[x=1, y=4,Lpos=90,Ldist=-3pt,Math,L=u_1]{A}
\Vertex[x=1, y=3,Lpos=0,Ldist=-1pt,Math,L=u_2]{B}
\Vertex[x=1, y=2,Lpos=-90,Ldist=-3pt,Math,L=u_3]{C}
\Vertex[x=2, y=4,Lpos=90,Ldist=-3pt,Math,L=v_1]{D}
\Vertex[x=2, y=3,Lpos=180,Ldist=-1pt,Math,L=v_2]{E}
\Vertex[x=2, y=2,Lpos=-90,Ldist=-3pt,Math,L=v_3]{F}
\Edges(B,C)
\Edges(E,F)
\Edges(B,F)
\Edges(C,F)
\Edges(C,E)
\Edge(D)(A)
\tikzset{EdgeStyle/.style = {-,bend right=90,color=gray}}
\Edge(A)(C)
\tikzset{EdgeStyle/.style = {-,bend left=90,color=gray}}
\Edge(D)(F)
\draw [dash pattern=on 3pt off 3pt] (A)-- (B);
\draw [dash pattern=on 3pt off 3pt] (D)-- (E);
\draw (2.4,3) node{$G_2$};
\draw (0.6,3) node{$G_1$};
\end{tikzpicture}}
\hskip 3cm
\subfloat[]{
\begin{tikzpicture}[x=2cm,y=2cm]
\renewcommand*{\VertexSmallMinSize}{1.2pt}
\renewcommand*{\VertexInnerSep}{1.2pt}
\GraphInit[vstyle=Classic]
\Vertex[x=1, y=4,Lpos=90,Ldist=-3pt,Math,L=u_1]{A}
\Vertex[x=1, y=3,Lpos=-15,Ldist=-3pt,Math,L=u_2]{B}
\Vertex[x=1, y=2,Lpos=-90,Ldist=-3pt,Math,L=u_3]{C}
\Vertex[x=2, y=4,Lpos=90,Ldist=-3pt,Math,L=v_1]{D}
\Vertex[x=2, y=2,Lpos=-90,Ldist=-3pt,Math,L=v_2]{E}
\Vertex[x=1.5, y=3,Lpos=160,Ldist=-3pt,Math,L=v_3]{F}
\Edges(B,C)
\Edges(E,F)
\Edges(B,F)
\Edges(C,F)
\Edges(C,E)
\Edge(D)(A)
\tikzset{EdgeStyle/.style = {-,bend right=90,color=gray}}
\Edge(A)(C)
\tikzset{EdgeStyle/.style = {-,bend right=30,color=gray}}
\Edge(D)(F)
\draw [dash pattern=on 3pt off 3pt] (A)-- (B);
\draw [dash pattern=on 3pt off 3pt] (D)-- (E);
\draw (1.85,3) node{$G_2$};
\draw (0.6,3) node{$G_1$};
\end{tikzpicture}}
\caption{The graph $G=g_4(G_1,G_2)$}
\label{fig4}
\end{figure}
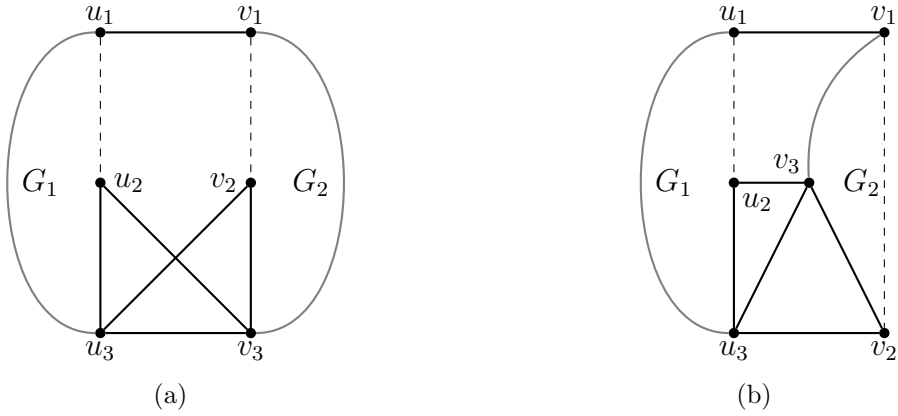

\medskip
Let $G_1, G_2$ be two planar graphs. Suppose that $u_1u_2u_3$ and $v_1v_2v_3$ are two paths each lie in the boundary of some face of $G_1$, $G_2$ respectively.  Let $G=g_4(G_1, G_2)$ be a graph (see Figure \ref{fig4}) obtained from $G_1$ and $G_2$ by

(a) Deleting the edges $u_1u_2$ and $v_1v_2$;

(b) Adding four new edges $u_1v_1$, $u_2v_3$, $u_3v_2$ and $u_3v_3$

\begin{figure}[H]
\centering
\subfloat[]{
\includegraphics[scale=1.]{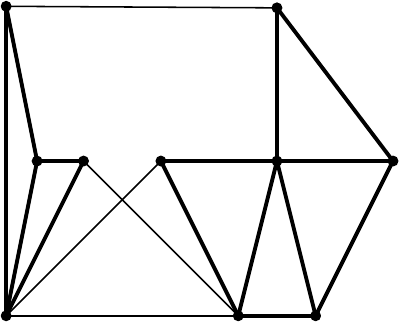}}
\hskip 3cm
\subfloat[]{
\includegraphics[scale=1.]{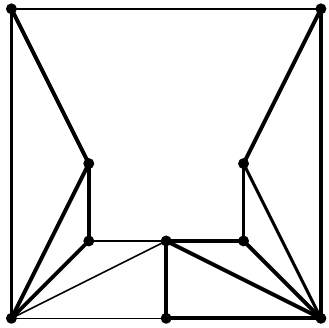}}
\caption{The graph $G=g_4(K_4, W_5)$}
\label{example4}
\end{figure}

If $G_1=K_4$ and $G_2=W_5$, the graph $G=g_4(K_4, W_5)$ is illustrated in Figure \ref{example4}. It is a 4-critical planar graph. Note that the graph in Figure \ref{example4} (b) is a redrawing of (a).

\begin{theorem}\label{thm2}
If $G_1$, $G_2$ and $G_{2k+1}$ are 4-critical planar graphs, then $G=g_4(G_1, G_2)$ is also a 4-critical planar graph.
\end{theorem}

{\bf Proof.} It is clear that $G=g_4(G_1, G_2)$ is planar, see Figure \ref{fig4}(b). Note that the graphs shown in Figure \ref{fig4}(a) and (b) are isomorphic, but their drawings are different. 

(i) We first show that $\chi (G)\ge 4$. Suppose on the contrary that $G$ is 3-colorable, then
since $G_1$ and $G_2$ are 4-critical, $G_1-u_1u_2$ and $G-v_1v_2$ are 3-colorable and $c(u_1)=c(u_2)$,$c(v_1)=c(v_2)=2$.  Without loss of generality, assume $c(u_1)=c(u_2)=1$ and $c(u_2)=c(u_3)=2$. Then $u_3$ and $v_3$ must be colored with color 3, this is impossible since $u_3$ and $v_3$ are adjacent. Therefore, we have $\chi (G)\ge 4$.

\medskip
(ii) We show that $\chi (G-e)\le 3$ for every edge in $G$. Consider the following  cases:

{\bf Case 1.} $e\in \{u_1v_1,u_3v_3, u_2v_3, u_3v_2\}$.

If $e=u_1v_1$, since $G_1$ is 4-critical, there is a 3-coloring $c_1$ of $G_1-u_1u_2$ such that $c_1(u_1)=c_1(u_2)=1$ and $c_1(u_3)=2$. Moreover, since $G_2$ is 4-critical, there is a 3-coloring $c_2$ of $G_2-v_1v_2$ such that $c_2(v_1)=c_2(v_2)=1$ and $c_2(v_3)=3$. It is clear that $c_1$ and $c_2$ extend to a 3-coloring of $G$.

If $e=u_3v_3$, since $G_1$ is 4-critical, there is a 3-coloring $c_1$ of $G_1-u_1u_2$ such that $c_1(u_1)=c_1(u_2)=1$ and $c_1(u_3)=2$. Moreover, since $G_2$ is 4-critical, there is a 3-coloring $c_2$ of $G_2-v_1v_2$ such that $c_2(v_1)=c_2(v_2)=1$ and $c_2(v_3)=2$. It is clear that $c_1$ and $c_2$ extend to a 3-coloring of $G$.

If $e=u_2v_3$ or $e=u_3v_2$, assume without loss of generality, that $e=u_2v_3$. Since $G_1$ is 4-critical, there is a 3-coloring $c_1$ of $G_1-u_1u_2$ such that $c_1(u_1)=c_1(u_2)=1$ and $c_1(u_3)=2$. Moreover, since $G_2$ is 4-critical, there is a 3-coloring $c_2$ of $G_2-v_1v_2$ such that $c_2(v_1)=c_2(v_2)=3$ and $c_2(v_3)=1$. It is clear that $c_1$ and $c_2$ extend to a 3-coloring of $G$.

{\bf Case 2.} $e\in E(G_1-u_1u_2)$ or $e\in E(G_2-v_1v_2)$.

Without loss of generality, assume that $e\in E(G_1-u_1u_2)$.

If $e=u_2u_3$, since $G_1$ is 4-critical, there is a 3-coloring $c_1$ of $G_1-e$ such that $c_1(u_2)=c_1(u_3)=1$ and $c_1(u_1)=2$. Moreover, since $G_2$ is 4-critical, there is a 3-coloring $c_2$ of $G_2-v_1v_2$ such that $c_2(v_1)=c_2(v_2)=1$ and $c_2(v_3)=2$. It is clear that $c_1$ and $c_2$ extend to a 3-coloring of $G$.

If $e\ne u_2u_3$, since $G_1$ is 4-critical, there is a 3-coloring $c_1$ of $G_1-e$ such that $c_1(u_1)=1$,$c_1(u_2)=2$ and $c_1(u_3)=k$ (where $k\in\{1,3\}$). Moreover, since $G_2$ is 4-critical, there is a 3-coloring $c_2$ of $G_2-v_1v_2$ such that $c_2(v_1)=c_2(v_2)=1$ and $c_2(v_3)=2$. It is straightforward to verify that $c_1$ and $c_2$ extend to a 3-coloring of $G$. \hfill$\Box$

\noindent{\bf Acknowledgments}

This research work is supported by National Natural Science Foundation of China under grant No. 11571168 and 11371193.


\begin{thebibliography}{10}
\bibitem{B.Zhou} H. L. Abbott, D. R. Hare and B. Zhou, Large faces in 4-critical planar graphs with minimum degree 4, 
{\em Combinatorica} 15 (4) (1995), 455-467.

\bibitem{AZ} H. L. Abbott and B. Zhou, The edge-density of 4-critical planar graphs, {\em Combinatorica} 11 (3) (1991), 185-189.

\bibitem{Bor}O. V. Borodin, Z. Dvo\v rákc, A. V. Kostochkad, B. Lidick\'y, M. Yancey, Planar 4-critical graphs with four triangles, {\em European Journal of Combinatorics} 41 (2014) 138-151.

\bibitem{hajos}
G. Haj\'{o}s, \"{U}ber eine Konstruktion nicht $n$-f\"{
a}rbbarer Graphen. Wiss. Z. Martin-Luther-Univ.Halle-Wittenberg, Math.-Naturw. Reihe 10, 116-117, 1961.

\bibitem{BG1} B. Gr\"unbaum, The edge-density of 4-critical planar graphs, {\em Combinatorica} 8 (1) (1988), 137-139.



\bibitem{TB}T. R. Jensen and B.Toft. Graph Coloring Problems. Wiley, New York, 1995.

\bibitem{Koester}G. Koester, 4-critical 4-valent planar graphs constructed with crowns, 
{\em Math. Scand.} 67 (1990), 15-22.

\bibitem{Koester1}G. Koester, On 4-critical planar graphs with high edge 
density, {\em Discrete Mathematics} 98 (1991), 147-151.

\bibitem{SS}H. Sachs and M. Stiebitz, On constructive methods in the theory of color-critical graphs, {\em Discrete Mathematics} 74 (1989), 201-226.









\end{thebibliography}
\end{document}